\documentclass[12pt,final]{amsart}

\usepackage{graphicx} 
\usepackage[margin=0.75in]{geometry}
\usepackage{amsmath}

\usepackage{verbatim}
\usepackage{amssymb}
\usepackage{amscd}

\newtheorem{lemma}{Lemma}[section]

\newtheorem{theorem}[lemma]{Theorem}

\theoremstyle{definition}

\theoremstyle{remark}
\newtheorem{remark}[lemma]{Remark}
\newtheorem{example}[lemma]{Example}

\newcommand{\N}{\ensuremath{{\mathbb N}}}

\newcommand{\R}{\ensuremath{\mathbb R}}

\author{Tatyana Barron}
\address{Department of Mathematics, University of Western Ontario, London Ontario N6A 5B7, Canada}
\email{tatyana.barron@uwo.ca}

\date{\today}

\thanks{}

\title{Geometric signals}
\begin{document}
\sloppy

\maketitle

\

\begin{abstract}
In signal processing, a signal is a function. Conceptually, replacing a function by its graph, and extending this approach to a more abstract setting, 
we define a signal as a submanifold $M$ of a Riemannian manifold 
(with corners) that satisfies additional conditions. In particular, it is a relative cobordism between two manifolds with boundaries. 
We define energy as the integral of the distance function to the first of these boundary manifolds. Composition of signals is composition of cobordisms. A "time variable" can appear explicitly if it is explictly given (for example, if the manifold is of the form $\Sigma\times [0,1]$). 
Otherwise, there is no 
designated "time dimension", although the cobordism may implicitly indicate the presence of dynamics. We interpret a local deformation of the metric as noise. The assumptions on $M$ allow to define a map $M\to M$ that we call a Fourier transform.      
We prove inequalities that illustrate the properties of energy of signals in this setting. 
 

\end{abstract}

\

\section{Introduction} 

This paper is about a geometric generalization of the concept of signal. The goal is to build an abstract mathematical model of  transmitting information via geometric objects which are, informally speaking, higher dimensional analogues of sound waveforms.     
We are not taking the most obvious path to defining a signal as an $\R^n$-valued function on an open subset of manifold. Instead, we 
are taking a more intuitive approach with cobordisms. This would allow to treat quite general geometric objects as information.

In practical applications, signal processing involves sampling, analog to digital conversion, quantization, and mathematical techniques that are needed because of  the hardware used to receive and analyze the signal.  We will concentrate our attention on the geometric concept of signal, rather than being concerned with digitizing such a signal or further processing the resulting data.  

\begin{figure}[hbp]
\includegraphics[width=2.5in]{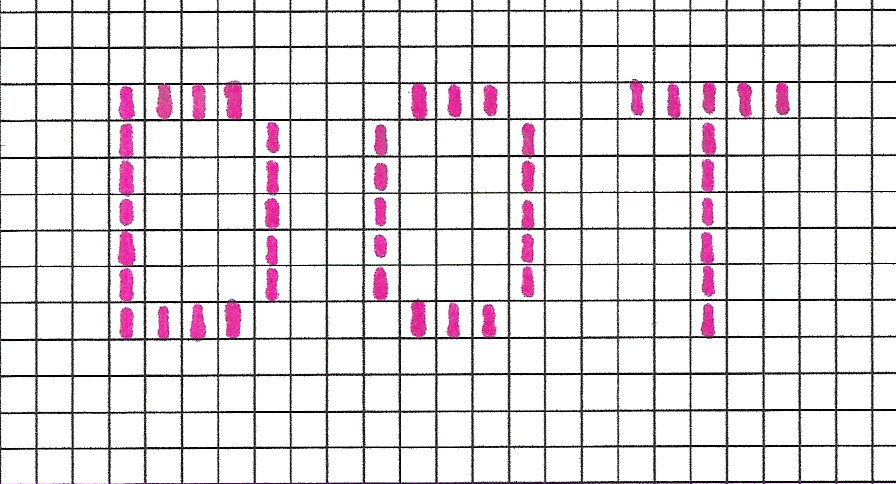}
\caption{The "empty squares" are the zero entries of the rectangular matrix. } \label{figdotpict}
\end{figure}

For background discussion, let's consider the word "dot". It can be transmitted as a finite sequence of three images, of the letters "D", "O", 'T", respectively. Each of the three images can be converted into digital data (e.g. by drawing the letters on the grid paper via shading appropriate squares and entering $1$ in the corresponding matrix for each shaded square, $0$ for every blank square Fig. \ref{figdotpict}). 

\begin{figure}[hbtp]
\includegraphics[width=2.5in]{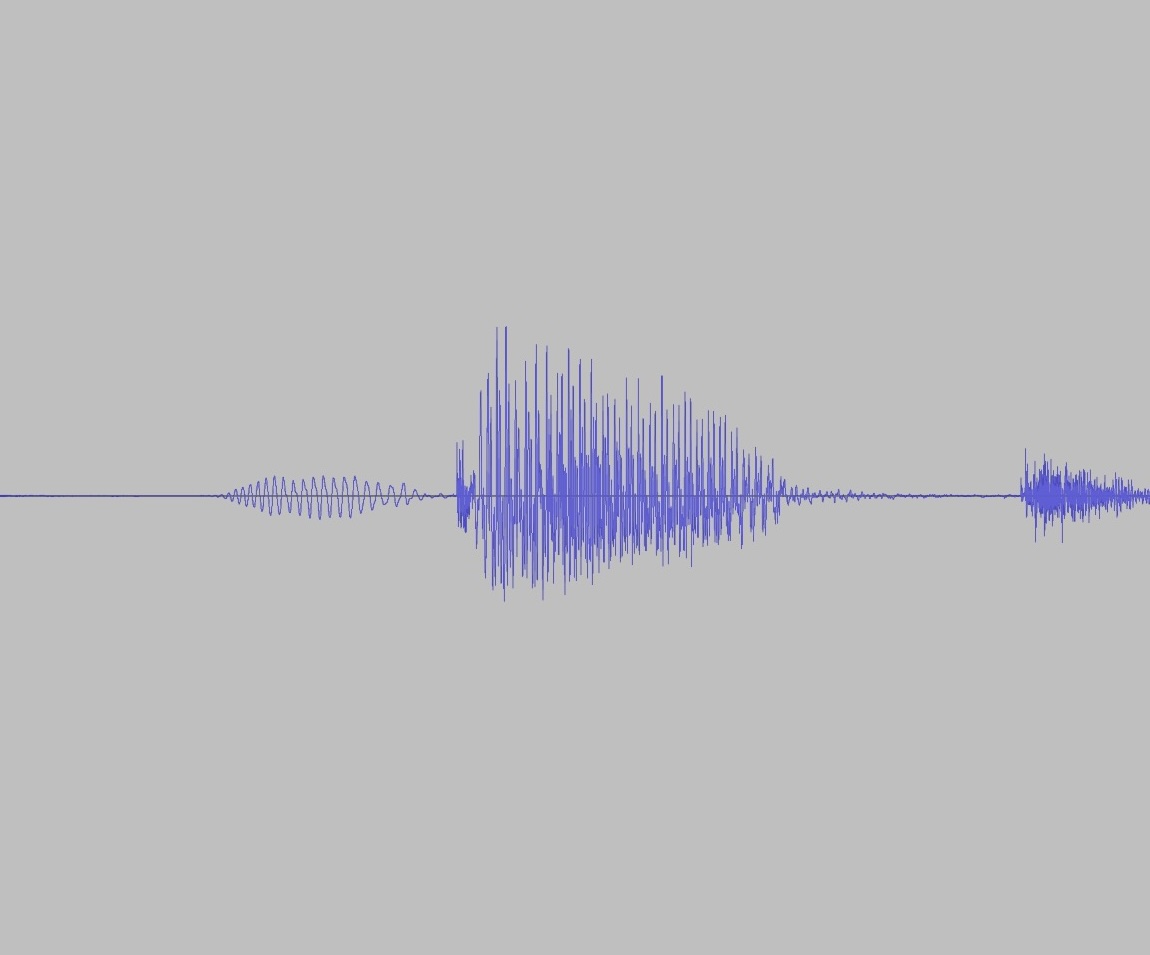}
\caption{A soundform of the word "dot". } \label{figdotso}
\end{figure}

Alternatively, the three letters can be represented by the binary representations of the numbers 4, 15, 20 (their positions in the English alphabet). Or, instead, one can send the word "dot" as an audio file ( a recording of a person saying this word), or the wave soundform of this recording (which can be subdivided into three parts, each for one of the three letters of the word Fig. \ref{figdotso}).   

 The mathematical aspects of these processes depend on the choice of sampling, quantization, encoding, transmission, other procedures related to speech processing, as well as the linguistics aspects such as the language and the alphabet (if the written language is based on an alphabet).       
For details see \cite{purdue2}, \cite{fbook}. 

Instead of doing all that, we can concentrate only on the geometric aspects and consider the immersed submanifolds of $\R^2$ (with the standard Riemannian metric)  that correspond to the data on Fig. \ref{figdotpict}, \ref{figdotso}. This is the general point of view in this paper. In order to account, intuitively,  for a "process" or "evolution" taking place, without explicitly defining a time variable, we use cobordisms.   

There is a vast amount of literature on applications of Riemannian geometry and its generalizations to signal processing, information theory, and computer science. It would be an impossible task to give a survey of this suject and to assign proper credit to all contributors. 
Here is an attempt to give a glimpse of this area. Work by Belkin and Niyogi, and co-authors, including \cite{bn:2003}, \cite{bn:2008}, \cite{bns:2006}, 
presents an intriguing vision. Broadly, one could view their approach as smoothing out the discrete data (or continuously approximating the discrete data). Typically the data sets are large and the data is put into a smooth manifold. Since every smooth manifold has a Riemannian metric, the analytic tools of Riemannian geometry are readily available. 
Similar reasoning supports the appearance of submanifolds in machine learning (see e.g. Ch. 6 \cite{machl}), with 
linear and nonlinear methods in dimensionality reduction. This also echoes the general philosophy of topological data analysis, most naively described as applying topological methods to discrete data mapped into a metric space.    
A particular kind of Riemannian manifolds is used in information geometry (see e.g. the discussion and references in Nielsen's survey \cite{nielsen}). In Menon's work \cite{menon} there are some interesting points about submanifolds, as well as regularity and numerical methods. In particular, he states that an embedding is an information transfer.    

Most generally, geometric aspects are often intrinsic to the numerical analysis, including effectiveness, accuracy, errors, in papers devoted to applications in a variety of areas, including, for instance, 
image processing, neural networks, or analysis of electromagnetic data  
(see \cite{karai}, \cite{ybi}, \cite{lchen}, \cite{wzhao} as some specific citations). 

In mathematics, the interplay between geometry (e.g. manifolds) and analysis (e.g. functions on this manifold)  is often pursued via proving  theorems 
that determine to what extent one determines the other (e.g. reconstruction theorems in noncommutative geometry, that allow to "rebuild" a manifold from its algebra of functions etc.; or no-go theorems, that say what can not happen on a manifold). An example of such work, with results that show how a manifold structure determines behaviour of certain functions on this manifold, is \cite{abarron}. 

In the present paper, we depart from this perspective. 
For example, instead of a function $f:\R\to\R$, we can consider its graph $\{ (t,f(t))|t\in\R\}$ which is a subset of the $ty$-plane. Generalizing this, 
we can consider an arbitrary curve in $\R^2$, or an arbitrary subset $C$ of the plane and consider this subset to be "data" or 'information" or "signal". There is no longer a need for it to be a graph of a function and there is no need to keep track of the global geometry of $\R^2$. 
Everything that we need to know about this curve should be localized on $C$. The setting in this paper involves cobordisms. 
In a follow up paper \cite{barronk} we use a more general definition, instead. We essentially 
view arbitrary geometric objects as signal/information and it also potentially allows for flexibility that, intuitively, should be needed to quantify local events that occur in neural networks.

\noindent {\bf Acknowldegements.} The valuable comments and suggestions from the reviewers are appreciated by the author. 
 
\section{Cobordisms as signals}

Let $k\ge 0$ be an integer. Let $X_1$ and $X_2$ be
closed $k$-dimensional oriented manifolds such that $X_1\cap X_2=\emptyset$ and $(X;X_1,X_2)$ is a $(k+1)$-dimensional (oriented) cobordism, i.e. $X$ is a compact $(k+1)$-dimensional oriented manifold 
with the boundary 
$$
\partial X=X_1 \sqcup X_2.
$$
\begin{figure}[hbtp]
\includegraphics[width=3.5in]{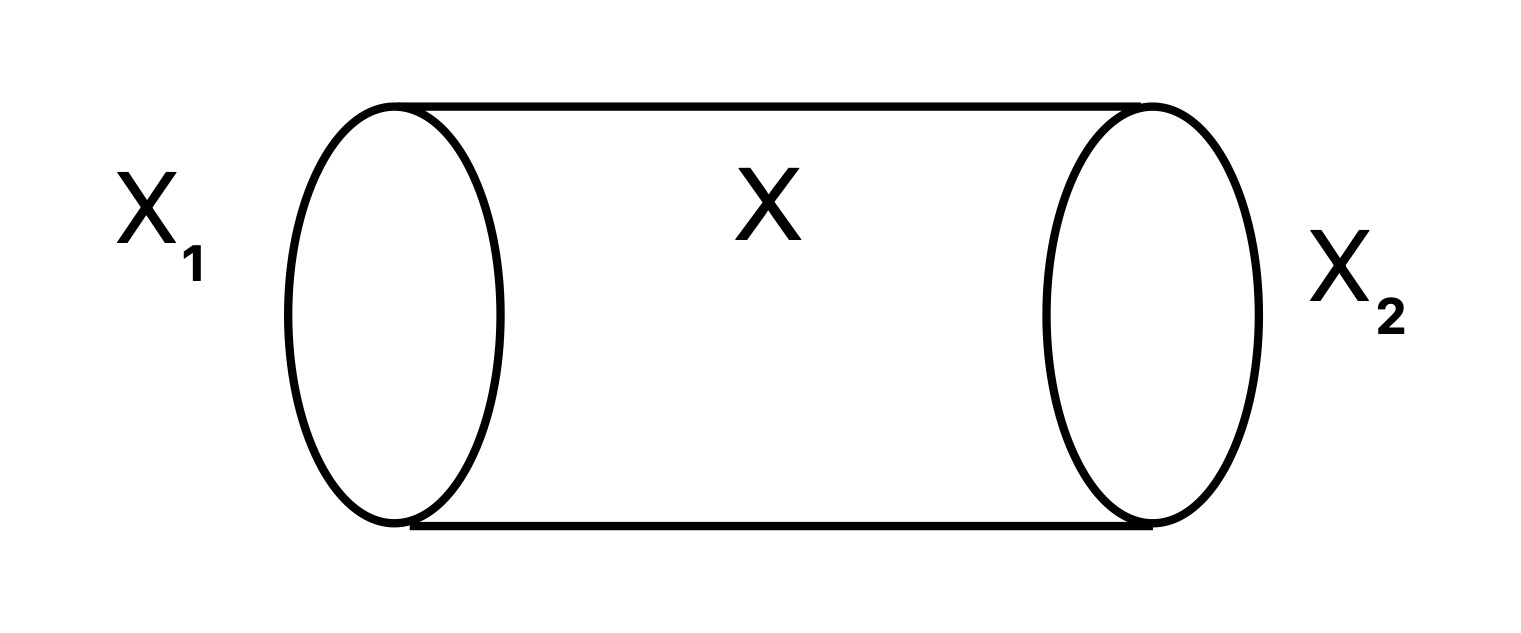}
\caption{A cobordism $(X;X_1,X_2)$.} \label{figx}
\end{figure} 
 Let $Y_1$ and $Y_2$ be
closed $k$-dimensional oriented manifolds and let $(Y;Y_1,Y_2)$ is a $(k+1)$-dimensional (oriented) cobordism. 
In general, all boundaries will be assumed to be nonempty where appropriate, and when we talk about a cobordism, as above, it will be assumed that 
$Y$, $Y_1$ and $Y_2$ are nonempty and $Y_1\cap Y_2=\emptyset$.
 
\begin{figure}[hbtp]
\includegraphics[width=3.5in]{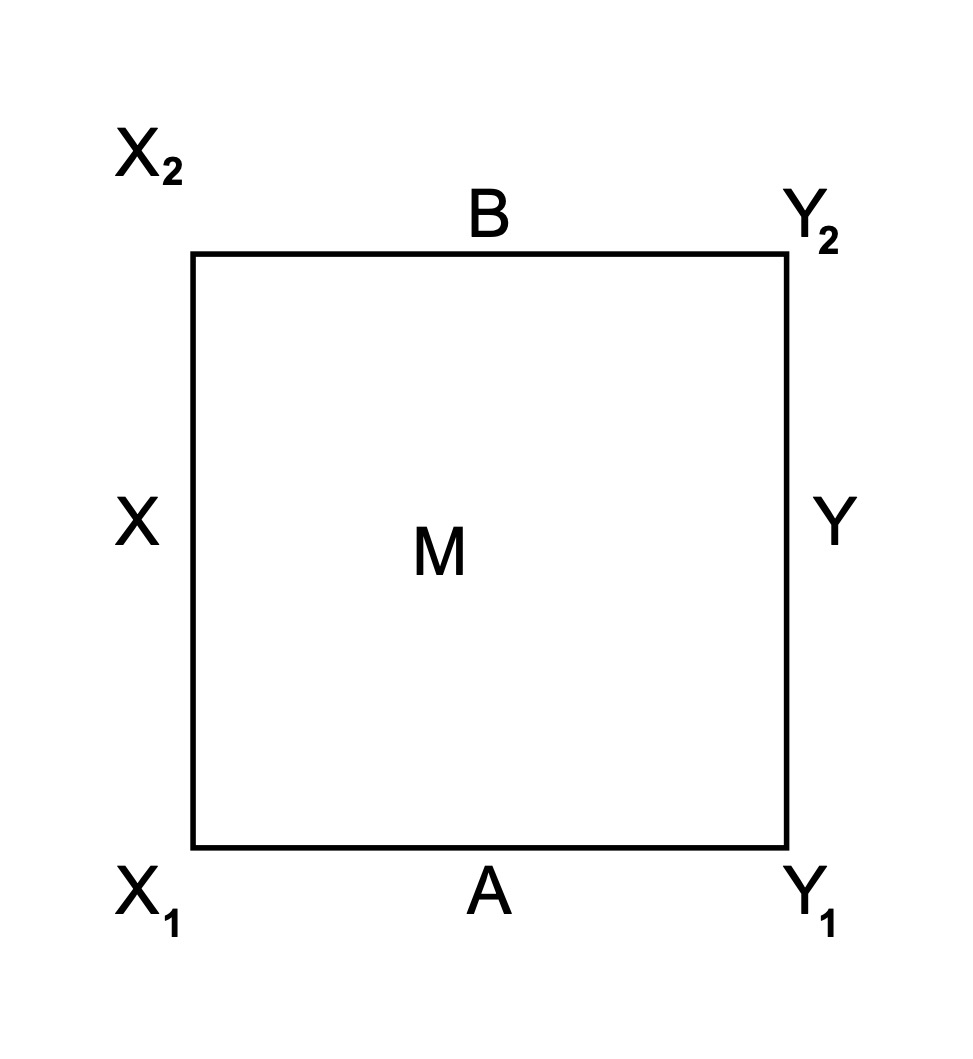}
\caption{A relative cobordism $(M;X,Y,\Sigma=A\cup B;\partial X, \partial Y)$, $k=0$.} \label{figsquar}
\end{figure}

\begin{figure}[hbtp]
\includegraphics[width=3.5in]{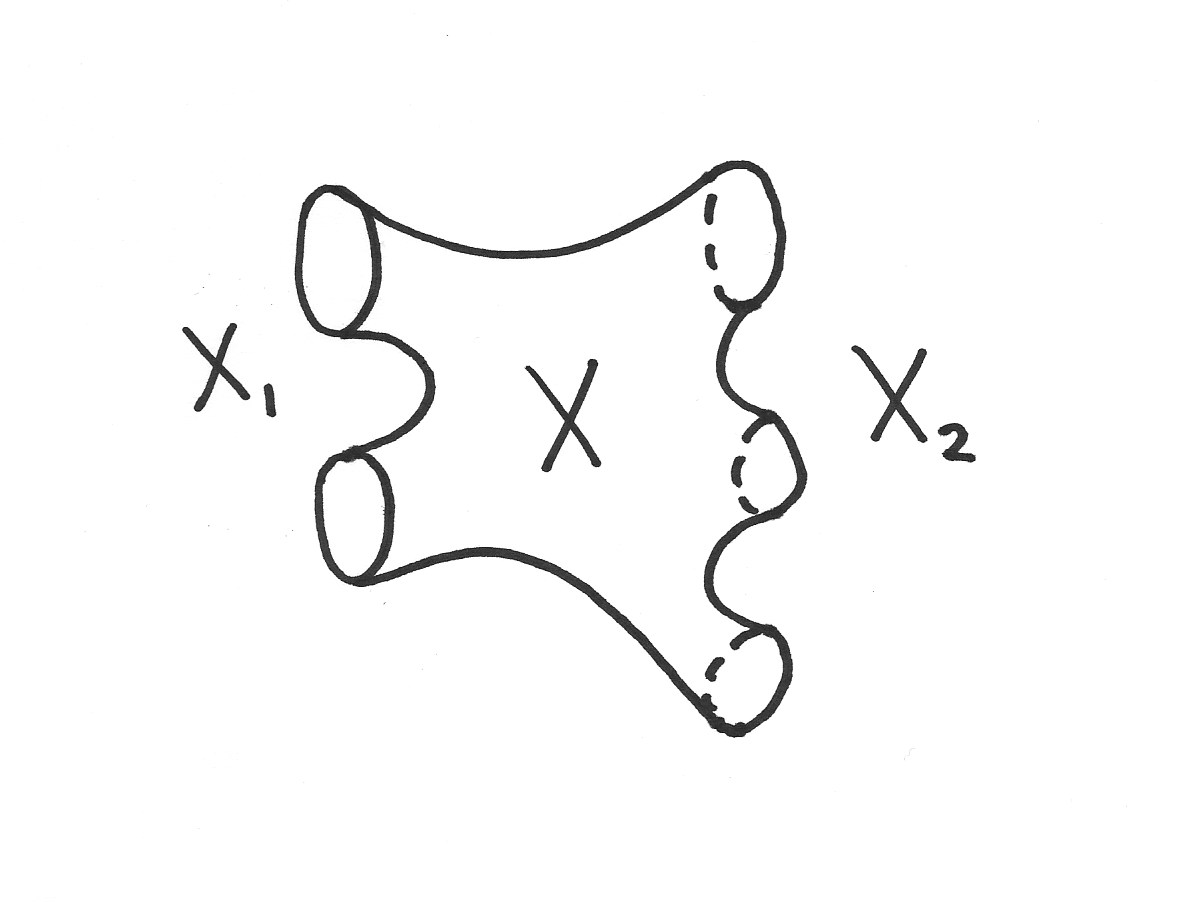}
\caption{A cobordism $(X;X_1,X_2)$, $n_1=2$, $n_2=3$.} \label{figxcirc}
\end{figure}

 Let $(M;X,Y,\Sigma;\partial X, \partial Y)$ be a $(k+2)$-dimensional cobordism 
between manifolds with nonempty boundaries (or, a relative cobordism) 
in the sense 
of \cite{bnr}, i.e. $(\Sigma, \partial X, \partial Y)$ is a $(k+1)$-dimensional cobordism, $\Sigma\ne \emptyset$, $X\cap Y= \emptyset$, 
$\Sigma\cap X= \partial X$, $\Sigma\cap Y= \partial Y$ and 
$$
\partial M= X \cup \Sigma \cup  Y.
$$ 
Assume, moreover, that
$$
\Sigma=A \sqcup B
$$
and  $(A;X_1,Y_1)$, $(B;X_2,Y_2)$ are $(k+1)$-dimensional cobordisms. 
$M$ is a manifold with corners \cite{laures}. Figure \ref{figsquar} shows an example where $k=0$, the boundary of $X$ consists of two points 
and the boundary of $Y$ consists of two points. 
\begin{example}
\label{exampleap}
For applications, a simplest typical example would be $k=1$, $X$ is a 
compact smooth surface with the boundary which is the disjoint union of $n_1+n_2$ circles ($n_1,n_2\in\N$), $X_1$ is a disjoint union 
of $n_1$ circles, $X_2$ is a disjoint union of $n_2$ circles, $Y=X$, $M=X\times [0,1]$. 
In Figure \ref{figx}, $n_1=n_2=1$.  
In Figure \ref{figxcirc}, $n_1=2$ and $n_2=3$.  Such $X$ could model a neuron or a part of a neural network. To take into account the electromagnetic field, one could consider submanifolds of $X\times [0,1]\times \R^6$. 
\end{example}

\begin{figure}[hbtp]
\includegraphics[width=3.5in]{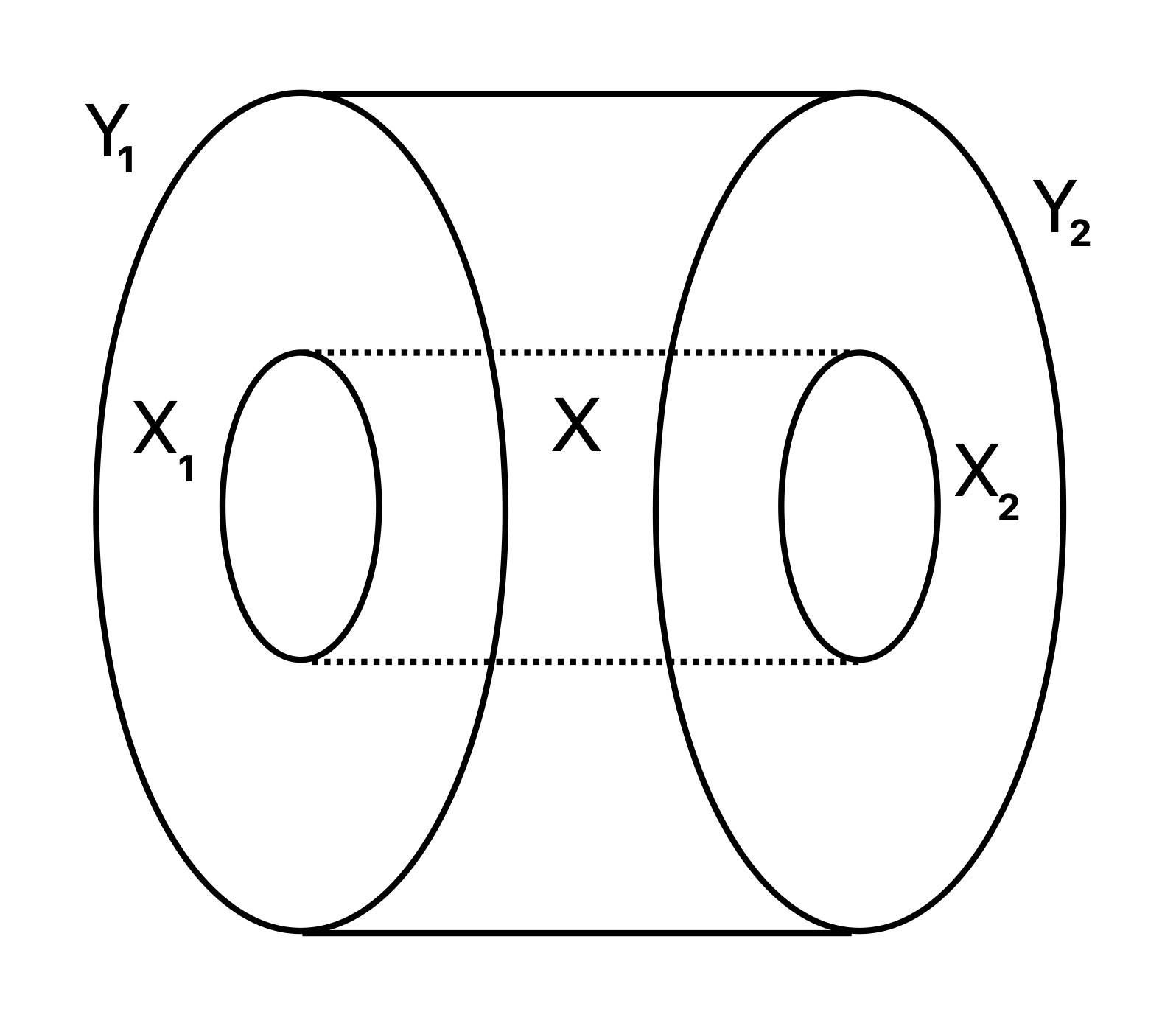}
\caption{A relative cobordism $(M;X,Y,\Sigma;\partial X, \partial Y)$.} \label{figxy}
\end{figure}

Let $(\tilde{M},\tilde{g})$ be a Riemannian manifold such that $M\subset \tilde{M}$ (a subset of $\tilde{M}$ which is an embedded submanifold via the inclusion map). We will write $M_g$ for $M$ equipped with the Riemannian metric $g$ induced by $\tilde g$. Unless explicitly stated otherwise, we will 
also assume that the Riemannian metric on every submanifold of $M$ is the one induced by $g$.    In Figure \ref{figsquar}, 
$M\subset \R^2=\tilde{M}$. 
Figure \ref{figxy} shows a relative cobordism with $X$ as in Example \ref{exampleap}, $k=1$, 
$n_1=n_2=1$, $M=X\times [0,1]$, and $M\subset\R^3=\tilde{M}$.
  
Denote by $dV_g$ the volume form of $g$ and by $\rho_g$ the Riemannian distance with respect to $g$. 
Define 
$$
f_X:M\to \R
$$
by $f_X(x)=\rho_g(x,X)$. Since $X$ is compact, $f_X$ is well defined and continuous on $M$. 
Similarly, define 
$$
f_A:M\to \R
$$
by $f_A(x)=\rho_g(x,A)$. 
\begin{remark}
Regularity of the distance function to a submanifold is treated in \cite{foote}. In \cite{foote}, the discussion is for $\R^n$ and it is noted that the proofs are similar for submanifolds of Riemannian manifolds. 
\end{remark} 
Denote $W=X\sqcup Y$.
Define 
\begin{itemize}
\item the  {\it energy of the signal} $M_g$ 
\begin{equation}
\label{energycob}
E(M_g)=\int_M f_A(x)dV_g(x)
\end{equation}
\item the {\it Fourier transform $F$ of the signal} $M_g$ 
as the map of cobordisms of manifolds with boundary 
$$
(M;X,Y,\Sigma;\partial X, \partial Y)\to (M;A,B,W;\partial X, \partial Y)
$$
induced by the identity map $M\to M$
\item the  {\it energy of the Fourier transform of the signal} $F(M_g)$ 
$$
E(F(M_g))=\int_M f_X(x)dV_g(x).
$$
\item a {\it noise} $(U,h)$ where $U\subset M$ is an open set and $h$ is a Riemannian metric on $M$ such that $h=g$ on $M-U$.  Then
$$
E(M_h)=\int_M\rho_h(x,A)dV_h(x)
$$ 
$$
E(F(M_h))=\int_M\rho_h(x,X)dV_h(x).
$$    

\item a {\it filter} a signal $M'=(M';X',Y',\Sigma';\partial X', \partial Y')$, such that $\partial X'=X_1' \sqcup X_2'$, $\partial Y'=Y_1' \sqcup Y_2'$, 
$\Sigma'=A' \sqcup B'$, 
$M'\subset M$, 
$X'\subset X$ and  $Y'\subset Y$.

\item composition of two signals $M=(M;X,Y,\Sigma=A\sqcup B;\partial X, \partial Y)$, $M'=(M';X',Y',\Sigma'=A' \sqcup B';\partial X', \partial Y')$ such that $Y=X'$ and $M\cap M'=X'$ is the signal
$$
M''=(M'';X,Y',\Sigma'';\partial X, \partial Y')
$$   
where 
$$
M''=M\cup M'
$$
$$
A''=A\cup A'
$$
$$
B''=B\cup B'
$$
$$
\Sigma''=\Sigma \cup \Sigma'=A''\sqcup B''
$$
We assume $A''\cap B''=\emptyset$. 
\end{itemize}
\begin{remark}
If there is no ambiguity about the metric and the metric is presumed to be the one induced by $g$, then 
we will sometimes omit $g$ from notation and write $E(M)=E(M_g)$, similarly for volume, diameter etc.   
\end{remark}
\begin{remark}
In signal processing,  
 a discrete signal is represented by a finite sequence of real numbers $(x_n)$. These come (via sampling)  from the sound waveform and should be understood as 
 the numbers that represent the air pressure at time $n$. The  signal is characterized by its
{\it energy}
$$
\sum_n|x_n|^2 
$$
and 
{\it magnitude} $\max \{ |x_n|\}$.
The energy of a continuous one-dimensional signal $y=x(t)$ is 
\begin{equation}
\label{energy}
\int |x(t)|^2 dt. 
\end{equation}
The concept of energy in Riemannian geometry is different and the value of energy for a curve in $\R^2$ is not the same integral as the integral (\ref{energy}) typically used in signal processing. One  can modify the signal $y=x(t)$ so that (\ref{energy}) becomes the Riemannian energy for the modified signal. 
 It is possible to write a similar argument that explains how our definition of energy (\ref{energycob}) relates to the other definitions, at least in a very simple low dimensional setting.
 We will discuss this in detail in upcoming work \cite{barronk}, where we will also consider the meaning of Fourier transform and noise (which is sometimes defined as the Fourier tranform of autocorrelation).       
\end{remark}
\begin{remark} 
To give an example of noise (noise as a local distortion of the metric, as stated in the list of definitions above), one can consider a small ball $B$ in $M$ and either a conformal deformation of the metric $g$ supported on $B$ or   a local diffeomorphism $\varphi$ that is the identity map on $M-Int(B)$, and the metric 
$h=\varphi^*g$.  
\end{remark} 
Recall (see e.g. section 3.2 \cite{hliu}) that the (boundary) injectivity radius $i_A>0$ of $A$ is defined as follows. 
Let $T^{\perp}A$ be the normal bundle of $A$ (in $M$), trivialized by the inward unit normal vector field $N$. Identify $T^{\perp}A$ with 
$A\times \R$ via this trivialization. For $p\in A$, denote by $\gamma$ the geodesic such that $\gamma(0)=p$ and $\gamma'(0)=N_p$. 
Let $D(p)=\inf \{ t>0|  \gamma_p(t)\in \partial M\}\in (0,\infty]$. The boundary exponential map 
$$
\exp ^{\perp}:(p,t)\mapsto \gamma_t(p)
$$ 
is defined on the subset of $A\times \R$ that consists of pairs $(p,t)$ with $0\le t<D(p)$. The number $\iota_A$ is defined as the supremum 
of $s\ge 0$ such that $\exp^{\perp}\Bigr | _{A\times [0,s)}$ is a diffeomorphism onto its image.   
The boundary injectivity radius $i_X>0$ of $X$ is defined similarly. 
\begin{theorem} 
\label{thmain1} Let 
$M=(M;X,Y,\Sigma=A\sqcup B;\partial X, \partial Y)$ be a signal. 

\noindent $(i)$  Then
$$
\frac{1}{ 1+\frac{ 4 \ {\mathrm{vol}}(M)  ( {\mathrm{diam}}(M)+{\mathrm{diam}}(A)+{\mathrm{diam}}(X))  }{ i_X^2 {\mathrm{vol}}(X)} }\le \frac{E(F(M))}{E(M)}\le 
$$
$$
1+\frac{ 4 \ {\mathrm{vol}}(M)  ( {\mathrm{diam}}(M)+{\mathrm{diam}}(A)+{\mathrm{diam}}(X))  }{ i_A^2 {\mathrm{vol}}(A)}. 
$$

\noindent $(ii)$  Let $(U, h)$ be noise, where  
$$
U=\{ x\in M| \ \rho_g(x,p)<\delta\} 
$$ 
for some $p\in M$ and  $\delta>0$. Let $0<\delta_0<\delta$ and let $0<\varepsilon<1$. 

Then there is a smooth function $a_{\varepsilon}:M\to\R$ such that 
\begin{equation}
\label{acond1}
a_{\varepsilon}(x)=\left\{ 
\begin{array}{ll }  
\varepsilon, \ {\mathrm{if}} \ x\in\{ x\in M| \ \rho_g(x,p)\le \delta_0\}  \\
1, \ {\mathrm{if}} \ x\in\{ x\in M| \ \rho_g(x,p)\ge \delta\} 
\end{array} \right. 
\end{equation}
\begin{equation}
\label{acond2}
0<a_{\varepsilon}(x)<1 \  {\mathrm{for \  all}}\ x\in\{ x\in M| \ \delta_0<\rho_g(x,p)< \delta\} , 
\end{equation}
and as $\varepsilon\to 0$ 
$$
\frac{E(F(M_{ha_{\varepsilon}})) }{E(M_{ha_{\varepsilon} })}=
\frac{\beta}{\gamma}\Bigl ( 1+C\varepsilon^{\frac{k+2}{2}} +O(\varepsilon^{k+2})\Bigr )
$$
where
$$
\beta=\int_{\{x\in M|\rho_g(p,x)>\delta_0\}}\rho_{ha_{\varepsilon}}(x,X)a_{\varepsilon}(x)^{\frac{k+2}{2}}dV_h(x)
$$
$$
\gamma=\int_{\{ x\in M|\rho_g(p,x)>\delta_0\}}\rho_{ha_{\varepsilon}}(x,A)a_{\varepsilon}(x)^{\frac{k+2}{2}}dV_h(x)
$$
$$
C=\frac{\int\limits_{\{x\in M|\rho_g(p,x)<\delta_0\}}\rho_{ha_{\varepsilon}}(x,X)dV_h(x) }{\beta}-
\frac{\int\limits_{\{x\in M|\rho_g(p,x)<\delta_0\}}\rho_{ha_{\varepsilon}}(x,A)dV_h(x) }{\gamma}.
$$
\noindent $(iii)$  Let $(U,h)$ be noise as in $(ii)$ and  
let $M'$ be a filter such that $X_1'=X_1$, $Y_1'=Y_1$, $A'=A$, $U\subset M-M'$. 
Then 
$$
E(M'_g)\le E(M_g)
$$
$$
E(M'_g)\le E(M_h).
$$
\end{theorem}
\begin{remark}
In Theorem \ref{thmain1} $(ii)$ the function $a_{\varepsilon}$, in a sense, modulates the noisy signal, from $M_h$ to $M_{ha_{\varepsilon}}$. 
In $(iii)$, $M'$ filters out the noise.  
\end{remark}
\begin{theorem}
\label{thmain2}
Let $M''$ be the composition of signals $M$ and $M'$. Then 
\begin{equation}
\label{ineq1}
E(M'')\le E(M)+E(M')
\end{equation}
\begin{equation}
\label{ineq2}
E(F(M''))\ge E(F(M)).
\end{equation}
\end{theorem}

\section{Proofs}

\paragraph{Proof of Theorem \ref{thmain1}.} 
By the mean value theorem, there is $a_0\in M$ such that 
$$
E(M)= \int_Mf_A(x)dV_g(x)=f_A(a_0){\mathrm{vol}}_g(M)
$$
and there is $x_0\in M$ such that 
$$
E(F(M))= \int_M f_X(x)dV_g(x)=f_X(x_0){\mathrm{vol}}_g(M).
$$

Since $A$ and $X$ are compact, there are $a_1\in A$ and $x_1\in X$ such that $\rho_g(a_0,a_1)=\rho_g(a_0,A)$ 
and $\rho_g(x_0,x_1)=\rho_g(x_0,X)$ . Let $b\in A\cap X$. Then  
$$
\frac{E(F(M))}{E(M)}= \frac{f_X(x_0)}{f_A(a_0)}=  \frac{\rho_g(x_0,X)}{f_A(a_0)}\le  \frac{\rho_g(x_0,b)+\rho_g(x_1,b)}{f_A(a_0)}\le
$$
$$
  \frac{\rho_g(a_0,x_0)+\rho_g(a_0,b)+\rho_g(x_1,b)}{f_A(a_0)}\le  
   \frac{\rho_g(a_0,x_0)+\rho_g(a_0,a_1)+\rho_g(a_1,b)+\rho_g(x_1,b)}{f_A(a_0)}\le
$$
\begin{equation}
\label{1plus}
1+\frac{ {\mathrm{diam}}(M)+{\mathrm{diam}}(A)+{\mathrm{diam}}(X) }{\rho_g(a_0,A)}.
\end{equation}
There is a neighborhood $A_1$ of $A$ in $M$  
and a diffeomorphism 
$$
\alpha:  A_1\to  A \times [0,i_A)
$$
defined by the normal exponential map. Let 
$$
A_2=\alpha^{-1}\Bigl ( A\times [\frac{i_A}{2},i_A)\Bigr ) .
$$
Let $\mathcal{U}=\{ U_1,...,U_m\}$ be an open cover of $A$ by manifold charts $\{ U_i, \varphi_i\}$ and let $\{ \psi_1,...,\psi_m\}$ be a smooth partition 
of unity subordinate to $\mathcal{U}$.  Then 
$$
f_A(a_0)=\frac{\int_Mf_A(x)dV_g(x)}{ {\mathrm{vol}}(M)}\ge \frac{1}{ {\mathrm{vol}}(M)} \int_{A_2} f_A(x)dV_g(x)=
$$
$$
\frac{1}{ {\mathrm{vol}}(M)} \int_{A_2} f_A(x)\sum_{j=1}^m\psi_j(x) dV_g(x)=
$$
\begin{equation}
\label{sumintegrals}
\frac{1}{ {\mathrm{vol}}(M)} \sum_{j=1}^m\int_{U_j\times [\frac{i_A}{2},i_A]} f_A(\alpha^{-1}(a,t))\psi_j(\alpha^{-1}(a,t) ) dV_g(\alpha^{-1}(a,t) ),
\end{equation}
where $a\in A$, $t\in [\frac{i_A}{2},i_A)$. 
Denote by $a_1^{(i)},...,a_{k+1}^{(i)}$ the coordinates in the $i$-th chart. Using the Fubini theorem, we get that (\ref{sumintegrals}) equals
$$
 \frac{1}{ {\mathrm{vol}}(M)}\int\limits_{ [\frac{i_A}{2},i_A]}  \sum_{i=1}^m
  \int\limits_{\varphi_i(U_i)} 
 f_A(\alpha^{-1}(\varphi_i^{-1}(a_1,...,a_{k+1}),t))\psi_i(\alpha^{-1}(\varphi_i^{-1}(a_1,...,a_{k+1}),t) ) 
 $$
 $$
 \sqrt{\det g(a_1,...,a_{k+1},t)}
  da_1 ...da_{k+1}dt.
 $$
 By the mean value theorem, there is $\frac{i_A}{2}<t_0<i_A$ such that 
 $$
 \int\limits_{  [ \frac{i_A}{2},i_A]}  
  \sum_{i=1}^m
 \int\limits_{\varphi_i(U_i)} 
 f_A(\alpha^{-1}(\varphi_i^{-1}(a_1,...,a_{k+1}),t))\psi_i(\alpha^{-1}(\varphi_i^{-1}(a_1,...,a_{k+1}),t) ) 
 $$
 $$
 \sqrt{\det g(a_1,...,a_{k+1},t)}
  da_1 ...da_{k+1}dt=
 $$
 \begin{equation}
 \label{ithint}
 \frac{i_A}{2}
  \sum_{i=1}^m
 \int\limits_{\varphi_i(U_i)} 
 f_A(\alpha^{-1}(\varphi_i^{-1}(a_1,...,a_{k+1}),t_0))\psi_i(\alpha^{-1}(\varphi_i^{-1}(a_1,...,a_{k+1}),t_0) ) 
\end{equation}
 $$
 \sqrt{\det g(a_1,...,a_{k+1},t_0)}
  da_1 ...da_{k+1}.
 $$
Then 
$$
f_A(a_0)\ge  \frac{1}{ {\mathrm{vol}}(M)}\frac{i_A}{2} \sum_{i=1}^m
 \int\limits_{\varphi_i(U_i)} 
 f_A(\alpha^{-1}(\varphi_i^{-1}(a_1,...,a_{k+1}),t_0))
 $$
 $$
 \psi_i(\alpha^{-1}(\varphi_i^{-1}(a_1,...,a_{k+1}),t_0) ) 
 \sqrt{\det g(a_1,...,a_{k+1},t_0)}
  da_1 ...da_{k+1}=
 $$
 $$
 \frac{1}{ {\mathrm{vol}}(M)}\frac{i_A}{2} \sum_{i=1}^m
 \int\limits_{U_i} 
 f_A(\alpha^{-1}(a,t_0))\psi_i(\alpha^{-1}(a,t_0) ) 
d\mu _A(a)=
 $$
$$
 \frac{1}{ {\mathrm{vol}}(M)} \frac{i_A}{2}\sum_{i=1}^m
 \int\limits_{A} 
 f_A(\alpha^{-1}(a,t_0))\psi_i(\alpha^{-1}(a,t_0) ) 
d\mu _A(a)=
$$
$$
 \frac{1}{ {\mathrm{vol}}(M)} \frac{i_A}{2}
 \int\limits_{A} 
 f_A(\alpha^{-1}(a,t_0))
d\mu _A(a)
 $$
 where $d\mu _A$ is the measure on $A$ induced by the Riemannian metric.  
 Hence 
 $$
\rho_g(a_0,A)=f_A(a_0)\ge  \frac{1 }{ {\mathrm{vol}}(M) }  \frac{i_A}{2}t_0{\mathrm{vol}}(A)\ge \frac{{\mathrm{vol}}(A) }{ {\mathrm{vol}}(M) }  \frac{i_A^2}{4}.
$$
Then, with (\ref{1plus}),  we get 
$$
\frac{E(F(M))}{E(M)}\le 1+\frac{ 4 \ {\mathrm{vol}}(M)  ( {\mathrm{diam}}(M)+{\mathrm{diam}}(A)+{\mathrm{diam}}(X))  }{ i_A^2 {\mathrm{vol}}(A)}.
$$
Repeating the argument for $\frac{E(M)}{E(F(M))}$, we get: 
$$
\frac{E(M)}{E(F(M))}\le 1+\frac{ 4 \ {\mathrm{vol}}(M)  ( {\mathrm{diam}}(M)+{\mathrm{diam}}(A)+{\mathrm{diam}}(X))  }{ i_X^2 {\mathrm{vol}}(X)}.
$$
This completes the proof of 
$(i)$.

Proof of $(ii)$. Let $\{ \varphi, \psi\}$ be a smooth partition of unity subordinate to the open cover 
$$
\{  \{x\in M|\rho_g(p,x)<\delta\}, \{x\in M|\rho_g(p,x)>\delta_0\} \} .
$$
Then the function $a_{\varepsilon}$ defined by 
$$
a_{\varepsilon}(x)=\varepsilon\varphi(x)+\psi(x)
$$ 
satisfies (\ref{acond1}), (\ref{acond2}). We have: 
$$
\frac{E(F(M_{ha_{\varepsilon}})) }{E(M_{ha_{\varepsilon} })}=
\frac{\int_M\rho_{ha_{\varepsilon}} (x,X)dV_{ha_{\varepsilon}} (x) } 
{\int_M \rho_{ ha_{\varepsilon}} (x,A)dV_{ha_{\varepsilon}} (x)}=
\frac{\int_M\rho_{ha_{\varepsilon}} (x,X) a_{\varepsilon}(x)^{\frac{k+2}{2}}dV_h (x) } 
{\int_M \rho_{ ha_{\varepsilon}} (x,A) a_{\varepsilon}(x)^{\frac{k+2}{2}}dV_h (x)}=
$$
$$
\frac{
\int_{\{ x\in M| \ \rho_g(x,p)> \delta_0\}}\rho_{ha_{\varepsilon}} (x,X) a_{\varepsilon}(x)^{\frac{k+2}{2}}dV_h (x) +\varepsilon^{\frac{k+2}{2}}
 \int_{\{ x\in M| \ \rho_g(x,p)\le \delta_0\}}\rho_{ha_{\varepsilon}} (x,X) dV_h (x) }
{\int_{\{ x\in M| \ \rho_g(x,p)> \delta_0\}}\rho_{ha_{\varepsilon}} (x,A) a_{\varepsilon}(x)^{\frac{k+2}{2}}dV_h (x) +\varepsilon^{\frac{k+2}{2}}
 \int_{\{ x\in M| \ \rho_g(x,p)\le \delta_0\}}\rho_{ha_{\varepsilon}} (x,A) dV_h (x) } =
 $$
 $$
 \frac{\beta}{\gamma}\frac{1+\varepsilon^{\frac{k+2}{2}} \frac{1}{\beta}
 \int_{\{x\in M|\rho_g(p,x)<\delta_0\}}\rho_{ha_{\varepsilon}}(x,X)dV_h(x) }
 {1+ \varepsilon^{\frac{k+2}{2}}\frac{1 }{\gamma}
 \int_{\{x\in M|\rho_g(p,x)<\delta_0\}}\rho_{ha_{\varepsilon}}(x,A)dV_h(x) }.
 $$
 The Maclaurin series for 
$$
\frac{1}{1+ \varepsilon^{\frac{k+2}{2}}\frac{1 }{\gamma}
 \int_{\{x\in M|\rho_g(p,x)<\delta_0\}}\rho_{ha_{\varepsilon}}(x,A)dV_h(x) }
$$
yields the desired statement.

 Proof of $(iii)$. 
 $$
E(M'_g)=\int_{M'}\rho_g(x,A')dV_g(x).
$$    
Since $A'=A$, $M'\subset M$, and $U\cap M'=\emptyset$, both inequalities follow. 
\qed

\paragraph{ Proof of Theorem \ref{thmain2}.} 
$$
E(M'')=\int_{M''}\rho_g(x,A'')dV_g(x)=\int_{M}\rho_g(x,A'')dV_g(x)+\int_{M'}\rho_g(x,A'')dV_g(x)
$$
Since for every $x$, $\rho_g(x,A'')\le \rho_g(x,A)$ and $\rho_g(x,A'')\le \rho_g(x,A')$, the inequality (\ref{ineq1}) follows. 
$$
E(F(M''))=\int_{M''}\rho_g(x,X)dV_g(x)=\int_{M}\rho_g(x,X)dV_g(x)+\int_{M'}\rho_g(x,X)dV_g(x)
$$
The inequality (\ref{ineq2}) follows. 
\qed 
\section{Conclusions}

In signal processing, a signal is  a function (typically, a function of time). See e.g. \cite{priem}, Chapter zero. 
In this paper, we define a signal to be a submanifold $M$ of Riemannian manifold, with extra conditions: it is a (relative) cobordism between two manifolds with boundary, $X$ and $Y$, and moreover, $M$ is also a cobordism between two other manifolds with boundaries, $A$ and $B$. We give examples. In particular, instead of a function of time, we would now consider its graph. This definition allows to define a self-map of $M$, which we call a Fourier transform. We define energy, which is a positive real number that characterizes a strength of the signal. Our definition is different from the standard Riemannian geometry definition and from the standard signal processing definition. We compare the three in the follow up paper \cite{barronk}. We say that noise is a local deformation of the metric, which makes sense since everything is bounded. 
We say that a modified signal filters out the noise if it is not affected by the noise. Composition of signals is composition of cobordisms. We prove energy inequalities. These inequalities give an idea about the behaviour of energy. This paper was intended as a step towards a geometric framework in the discussion of signals. We discuss applications in \cite{barronk}.

%

%
%

\end{document}